\documentclass{article}
\usepackage[]{epsf,epsfig,amsmath,amssymb,amsfonts,latexsym}
\title{An efficiency upper bound for inverse covariance estimation}
\author{Ronen Eldan \thanks{Supported in part by the Israel Science Foundation and by a Marie Curie Grant from the Commission of the European Communities}}
\newtheorem{theorem}{Theorem}[section]

\newtheorem{remark}{Remark}[section]

\def \PP {\mathbb P}
\def\qed{\hfill $\vcenter{\hrule height .3mm
\hbox {\vrule width .3mm height 2.1mm \kern 2mm \vrule width .3mm
height 2.1mm} \hrule height .3mm}$ \bigskip}
\def\P{\mathbb{P}}

\def\EE{\mathbb{E}}

\def\RR{\mathbb{R}}

\def\Sph{S^{d-1}}

\begin{document}
\maketitle

\begin{abstract}
We derive a quantitative upper bound for the efficiency of estimating entries in the inverse covariance matrix of a high dimensional distribution.
We show that in order to approximate an off-diagonal entry of the density matrix of a $d$-dimensional Gaussian random vector, one needs at least a number of samples proportional to $d$. Furthermore, we show that with $n \ll d$ samples, the hypothesis that two given coordinates 
are fully correlated, when all other coordinates are conditioned to be zero, cannot be told apart from the hypothesis that the two are uncorrelated.
\end{abstract}
\bigskip
\section{Introduction}
The problem of estimating the population covariance matrix given a sample of $n$ i.i.d. observations $X_1,...,X_n$ in $\RR^d$
has been extensively studied. Estimation of covariance matrices plays a key role in many data analysis techniques (e.g. in principal component analysis, discriminant analysis, graphical models). \\ \\
It has been shown in \cite{ALPT} that when the measure is log-concave, the empirical covariance matrix gives a good approximation when $n=\Omega(d)$. In the case $n<d$, it is clear that the empirical covariance matrix cannot give a good approximation for the population covariance matrix, since it is not of full rank. However, a-priori, we could hope that other approximation schemes may still work. Later in this note, we will see that it is not the case. \\ \\
An easier goal than approximating the entire convariance matrix $A$ would be to approximate a single entry in $A^{-1}$. 
The latter has a rather natural interpretation: 
Given a multivariate Gaussian random vector $Y=(Y_1,...,Y_d)$ whose covariance matrix is $A$ (namely, for all $1 \leq i,j \leq d$, one has $\EE[Y_i Y_j] = A_{i,j}$), and two indices $1 \leq k_1,k_2 \leq d$, define
\begin{equation} \label{alphas}
\alpha_{i,j} = \lim_{\epsilon \to 0} \EE[Y_i Y_j ~~ | ~~ |Y_k| < \epsilon, \forall k \notin \{k_1, k_2 \} ],
\end{equation}
for all $i,j \in \{k_1, k_2\}$. One may interpret the quantity $\alpha_{i,j}$ as the effective correlation between $Y_i$ and $Y_j$, 
in the sense that it neutralizes "indirect" effects caused by correlation with a third variable $Y_k$, $k \notin \{k_1,k_2\}$.
Now, it is easily seen that when $A$ is invertible, there is a simple relation the numbers $\alpha_{i,j}$ and 
the matrix $A^{-1}$, namely,
$$
\left (
\begin{array}{cc}
\alpha_{i,i} ~ \alpha_{i,j} \\
\alpha_{j,i} ~ \alpha_{j,j}
\end{array}
\right )^{-1} = 
\left (
\begin{array}{cc}
(A^{-1})_{i,i} ~ (A^{-1})_{i,j} \\
(A^{-1})_{j,i} ~ (A^{-1})_{j,j}
\end{array}
\right ).
$$

As an example, if the indices represent a set of genes, and the quantity $Y_i$ represents presence or absence of the $i^{th}$ gene, biologists are often interested to know whether or not a certain correlation between the presence of two different genes is due to the fact that both genes depend on a third gene. The number $\alpha_{i,j}$ gives an indication to whether these two genes are directly correlated, rather than being both correlated with a third gene. \\ \\
The goal of this short note is to introduce an {\it information-theoretic} lower bound for the above question, and show that the number of samples
needed in order to estimate the numbers $\alpha_{i,j}$ above is essentially the same as the minimum number of samples needed to estimate the entire population covariance matrix using the empirical covariance matrix. \\ \\
Before we formulate the result, let us introduce some notation. Fix a dimension $d$, and consider the Euclidean space $\RR^d$, and its standard basis $e_1,...,e_d$. Define $E = span \{e_1, e_2 \}$ and let $P_E$ be the orthogonal projection onto $E$.
Let $X$ be a standard Gaussian random vector in $\RR^d$. Note that for a positive semi-definite symmetric matrix $A$, the covariance matrix of the random vector $A^{1/2} X$ is exactly $A$. \\ \\
Next, denote by $B_d$ the Euclidean unit ball in $\RR^d$, and for a symmetric matrix $A \in GL(d)$, define $C_E(A)$ to be the covariance matrix of the uniform distribution on the ellipse $A^{1/2} B_d \cap E$ (whose dimension is between 0 and 2). Observe that the density of the vector $A^{1/2} X$ is constant on ellipsoids of the form $t A^{1/2} \partial B_d$, $t > 0$ and that there exists a constant $K_d$ depending only on the dimension $d$ such that
$$
K_d C_E(A) = \left (
\begin{array}{cc}
\alpha_{1,1} ~ \alpha_{1,2} \\
\alpha_{2,1} ~ \alpha_{2,2}
\end{array}
\right ),
$$
where $\alpha_{i,j}$ are the constants defined in equation \ref{alphas} (with $Y = A^{1/2} X$ and $k_1=1,k_2=2$). Our main goal boils down to showing that the entries of the matrix $C_E(A)$ cannot be approximated with a reasonable probability. We prove the following theorem:
\begin{theorem}
Suppose $n<\frac{d}{3}$. There does not exist a function $F:(\RR^d)^n \to \{0,1,2 \}$ such that for every positive semi-definite matrix $A \in GL(d)$, one has
\begin{equation} \label{finalgoal}
\P \left (\left \{ F(A^{1/2} X_1, ..., A^{1/2} X_n) = rank(C_E(A)) \right \} \right ) > 0.9
\end{equation}
where $X_1,...,X_n$ are independent standard Gaussian random vectors in $\RR^d$.
\end{theorem}
In other words, given $\frac{d}{3}$ samples or less, not only we cannot approximate the constants $\alpha_{i,j}$, but we cannot even determine the rank of the matrix $C_E(A)$ with a reasonable probability. \\ \\
The idea of the proof is the following: Let $X_1,...,X_n$ be independent standard Gaussian random vectors. We construct two random covariance matrices $A,B$ such that almost surely, $rank C_E(A) \neq rank C_E(B)$. On the other hand, the random matrices $A$ and $B$ will be constructed such that the total variation distance between the two following distributions on $(\RR^d)^n$ will be rather small: the first distribution is attained by randomly generating an instance of $A$ and then considering the sequence $(A^{1/2} X_1,..., A^{1/2} X_n)$, and the second by doing the same, replacing $A$ with $B$. Note that, conditioning on $A$ and $B$, the above are sequences of independent samples. A small total variation distance implies that for every function $F$, the total variation distance between the random variables $F(X_1,...,X_n)$ and $F(Y_1,...,Y_n)$ will be rather small, which means that no function $F$ can distinguish between the two. \\ \\

It is interesting to inspect the result of this note in view of some positive results concerning the estimation of the covariance matrix which appeared recently. The results provide methods to approximate, or partly approximate the covariance matrix or its inverse when some extra assumptions about the distribution of $X$ can be made. For example, when the covariance matrix is assumed to be rather sparse, some methods can be used in order to estimate a symmetric part of it, as in \cite{LV}, or the inverse matrix, as in \cite{BLRZ}, given a rather small number of samples. See also
\cite{V} for background and more related results. \\ \\

\begin{it} Acknowledgements \end{it} The author would like to thank Bo`az Klartag and Roman Vershynin for introducing him to the question and for fruitful discussions, and would also like to thank the anonymous referee for several useful comments, corrections and suggestions.

\section{Proof of the theorem}

To prove the theorem, we assume by contradiction that there exists a function $F: (\RR^d)^n \to \{0,1,2 \}$ satisfying (\ref{finalgoal}). \\ \\
We begin with the construction of two families of Gaussian vectors: \\
Let $X_1,...,X_n, \tilde Y_1,...,\tilde Y_n$ be independent samples of the standard Gaussian vector in $\RR^d$, and let 
let $\theta$ be a random variable uniformly distributed on $\Sph$ and independent from the above. Define $A = Proj_{\theta^\perp}$ and
$Y_i = A^{1/2} \tilde Y_i$ for $1 \leq i \leq n$. Clearly, when conditioned on $\theta$, $Y_1,...,Y_n$ are independent samples of some Gaussian distribution. It follows from the definition of $A$ that $C_E(A)$ is of rank 1 whenever $\theta \notin E^\perp$, which means that,
$$
\PP(rank(C_E(A)) = 1 ) = 1
$$ 
and by the assumption (\ref{finalgoal}) along with the conditional independence of $Y_1,...,Y_n$ with respect to $\theta$, it follows that
\begin{equation} \label{FYgood}
\PP(F(Y_1,...,Y_n) = 1) = \EE_\theta ( \PP(F(Y_1,...,Y_n) = 1) | ~ \theta ) > 0.9.
\end{equation}
Our next step is to show that if $F$ satisfies the assumption (\ref{finalgoal}) then there must also exist a function $G:(\RR^d)^n \to \{0,1,2\}$, invariant under the action of $SO(d)$, which satisfies a slightly weaker version of (\ref{finalgoal}). 
To this end, let $T$ be a random orthogonal matrix distributed uniformly according to the Haar measure on $SO(d)$ and independent of all the above. By the construction of the sequences, it is clear that 
$$
(T(X_1),...,T(X_n)) \sim (X_1,...,X_n) \mbox{ and } (T(Y_1),...,T(Y_n)) \sim (Y_1,...,Y_n).
$$
The assumption (\ref{finalgoal}) and equation (\ref{FYgood}) now give
\begin{equation}
\PP(F(T(X_1),...,T(X_n)) = 2) > 0.9, ~~ \PP(F(T(Y_1),...,T(Y_n)) = 1) > 0.9.
\end{equation}
Therefore, denoting
$$
G(Z_1,...,Z_n) = 
\begin{cases}
2 &,\; \EE_T (F(T(Z_1),...,T(Z_n))) > \frac 3 2 \\
1 &,\; \frac{1}{2} \leq \EE_T (F(T(Z_1),...,T(Z_n))) < \frac 3 2 \\
0 &,\; \mbox{otherwise}\\
\end{cases}
$$
it is easily checked that $G$ will satisfy:
\begin{equation} \label{goodalg}
\PP(G(X_1,...,X_n) = 2) > 0.8, ~~ \PP(G(Y_1,...,Y_n) = 1) > 0.8.
\end{equation}
The \emph{total variation distance} between two random variables $X,Y$ with values in $W$ is defined as
$$
d_{TV}(X,Y) = \sup_{A \subset W} | \PP(X \in A) - \PP(Y \in A) |.
$$
Equation (\ref{goodalg}) implies that,
\begin{equation} \label{totalvar}
d_{TV} (G(X_1,...,X_n), G(Y_1,...,Y_n)) > 0.6.
\end{equation}
Since $G$ is invariant under rotations, and since one can always choose an orthogonal transformation $T$ such that
$$
T(X_i) \in span \{e_1,...,e_i \}, ~~ \forall 1 \leq i \leq n,
$$
it is clear that the function $G$ must only depend on the Gram matrix of the samples. So,
$$
d_{TV} (G(X_1,...,X_n), G(Y_1,...,Y_n)) \leq d_{TV}(Gr(X_1,...,X_n),Gr(Y_1,...,Y_n))
$$
where $Gr(\cdot)$ denotes the Gram matrix. \\
Clearly,
$$
Gr(X_1,...,X_n) \sim W_n(Id, d)
$$
where $W_n(C, p)$ is the Wishart distribution of dimension $n$ with $p$ degrees of freedom and covariance matrix $C$.  \\ \\
Next, let us try to understand the distribution of $Gr(Y_1,...,Y_n)$. To that end, we make the following observation: let $\varphi_1,...,\varphi_{d-1}$ be random vectors in $\RR^d$ such that $\{ \varphi_1,..,\varphi_{d-1}, \theta \}$ is almost surely an orthonormal basis of $\RR^d$. By the construction of the vectors $Y_1,...,Y_n$, one may write
$$
Y_i = \sum_{k=1}^{d-1} \Gamma_{i,k} \varphi_k, ~~ \forall 1 \leq i \leq n
$$
where $\{ \Gamma_{i,k} \}_{i,k=1}^\infty$ is an infinite matrix of independent standard Gaussian variables, independent from $\{ \varphi_1,..,\varphi_{d-1}, \theta \}$. So one has, for all $1 \leq i,j \leq n$,
$$
\langle Y_i, Y_j \rangle = \sum_{k=1}^{d-1} \langle Y_i, \varphi_k \rangle \langle Y_j, \varphi_k \rangle = \sum_{k=1}^{d-1} \Gamma_{i,k} \Gamma_{j,k}.
$$
By definition of the Wishart matrix, we get
$$
Gr(Y_1,...,Y_n) \sim W_n(Id, d-1).
$$
Our task is therefore to estimate,
$$
d_{TV}(W_n(Id, d-1), W_n(Id, d)).
$$
$$
~~
$$
Let $\mathcal{P} \subset \RR^{n^2}$ be the cone of positive semi-definite matrices. It is well known (see e.g., \cite{W}) that a random matrix $A \sim W_n(Id, d)$ has the following density with respect to the Lebesgue measure on $\mathcal{P}$:
$$
f_{n,d}(A) := \frac{\det(A)^{\frac 1 2 (d-n-1)} \exp(- \frac 1 2 Trace(A))}{2^{\frac 1 2 dn} \pi^{\frac 1 4 n (n-1)} \prod_{i=1}^n \Gamma(\frac 1 2 (d+1-i)) }
$$
whenever $n \leq d$. Denote the measure expressing the law of $W_n(Id, d)$ by $\mu_{n,d}$.
We would like to estimate the total variation metric between $\mu_{n,d}$ and
$\mu_{n,d-1}$. For this, we write,
$$
d_{TV}(W_n(Id, d-1), W_n(Id, d)) = \frac{1}{2} \int_{\RR^{n^2}} |f_{n,d-1}(A) - f_{n,d}(A)| d \lambda(A)
$$
where $\lambda$ is the lebesgue measure on $\mathcal{P}$. Denote,
$$Z(n,d) = 2^{\frac 1 2 dn} \pi^{\frac 1 4 n (n-1)} \prod_{i=1}^n \Gamma \left (\frac 1 2 (d+1-i) \right ),$$
so that
$$
\int_{\mathcal{P}} \det(A)^{\frac 1 2 (d-n-1)} \exp \left (- \frac 1 2 Trace(A) \right ) = Z(n,d).
$$
We have,
$$
d_{TV}(W_n(Id, d-1), W_n(Id, d)) = \frac{1}{2} \int_{\mathcal{P}} \left |1 - \det(A)^{1/2} \frac{Z(n,d-1)}{Z(n,d)}  \right| d \mu_{n,d-1} (A).
$$
Note that,
$$
\int_{\mathcal{P}} \det(A)^{1/2} d \mu_{n,d-1} (A) = \frac{Z(n,d)}{Z(n,d-1)}.
$$
The last two equations give,
$$
d_{TV}(W_n(Id, d-1), W_n(Id, d)) = 
$$
$$
\frac{1}{2} \int_{\mathcal{P}} \left |1 - \frac{\det(A)^{1/2}}{\int_{\mathcal{P}} \det(A)^{1/2} d \mu_{n,d-1} (A) } \right| d \mu_{n,d-1} (A) \leq 
$$
$$
\frac{1}{2} \sqrt { \int_{\mathcal{P}} \left (1 - \frac{\det(A)^{1/2}}{\int_{\mathcal{P}} \det(A)^{1/2} d \mu_{n,d-1}} \right)^2 d \mu_{n,d-1}(A)} =
$$
$$
\frac{1}{2} \frac{\sqrt{Var[\det^{1/2}(W_n(Id, d-1))]}}{\EE[\det^{1/2}(W_n(Id, d-1))]}.
$$
Let $X$ be a random variable such that $\EE[|X|^4]$ exists. It follows from Lyapunov's inequality (see e.g., \cite[P. 117]{PPT}) that
$$
\frac{\EE[|X|^4]}{\EE[|X|]} \geq \left ( \frac{\EE[|X|^2]}{\EE[|X|]} \right )^3
$$
So,
$$
\EE \left [|X|^4 \right ] - \EE \left [|X|^2 \right ]^2 \geq \frac{\EE \left [|X|^2 \right ]^3}{\EE[|X|]^2} - \EE \left [|X|^2 \right ]^2 
$$
And so,
$$
\frac{ Var \left [|X|^2 \right ] }{ \EE \left [|X|^2 \right ]^2 }  \geq \frac{ Var[|X|] }{ \EE \left [|X| \right ]^2 }.
$$
Using this inequality with the random variable $X \sim \det^{1/2}(W_n(Id, d-1))$ gives,
\begin{equation}
d_{TV}(W_n(Id, d-1), W_n(Id, d)) \leq \frac{1}{2}  \frac{\sqrt{ Var[\det(W_n(Id, d-1))]}}{\EE[\det(W_n(Id, d-1))]} .
\end{equation}
As shown in \cite[Theorem 4.4]{DMO}, one has
$$
Var[\det(W_n(Id, d-1))] = \frac{(d-1)!}{(d-1-n)!} \left (\frac{((d-1)+2)!}{((d-1)+2-n)!} - \frac{(d-1)!}{((d-1)-n)!} \right )
$$
and,
$$
\EE[\det(W_n(Id, d-1)] = \frac{(d-1)!}{(d-1-n)!}.
$$
So,
$$
d_{TV}(W_n(Id, d-1), W_n(Id, d)) \leq \frac{1}{2} \sqrt {\frac{(d+1)! (d-1-n)!}{ (d+1-n)! (d-1)! } - 1 } =
$$
$$
\frac{1}{2} \sqrt {\frac{(d)(d+1)}{(d-n)(d-n+1)} - 1}
$$
The assumption $n<\frac{d}{3}$ implies that the above expression is smaller than $0.6$.
This contradicts (\ref{totalvar}) and the proof is finished. \qed

\begin{remark}
It is not hard to see that if $n < \frac{d}{3}$ then there actually exists a (deterministic) matrix $A \in GL(d)$ such that if $Z_1,...,Z_n$ are independent standard Gaussian random vectors in $\RR^d$, one has,
$$
\PP(F(A^{1/2} Z_1,...,A^{1/2} Z_n) = rank C_E(A) ) < 0.9.
$$
Indeed, define $A_\theta = Proj_{\theta^\perp}$ and let $\Theta$ be a uniform point on $\Sph$. We have shown that
$$
\EE_\Theta[ \PP(F(A_\Theta^{1/2} Z_1,...,A_\Theta^{1/2} Z_n) = rank C_E(A_\Theta)) ] < 0.9.
$$
This implies that there exists a specific choice of $\theta \in \Sph \setminus E^\perp$ such that,
$$
\PP(F(A_\theta^{1/2} Z_1,...,A_\theta^{1/2} Z_n) = rank C_E(A_\theta)) < 0.9.
$$
The two Gaussian random vectors $Z, A_\theta^{1/2} Z$ are thus, in some sense, indistinguishable by $F$.
\end{remark}

\begin{remark}
It can be checked that when $n \ll d$, the function $F$ cannot do much better than being correct with probability 
$\frac{1}{3}$, hence, it cannot do better than guessing the rank of $C_E(A)$.
\end{remark}

\begin{remark}
Following the same lines of proof, one can also show that the correlation between two coordinates cannot be approximated also 
when conditioning on all but $k$ coordinates to be zero, where $k$ is a constant and $d \to \infty$.
\end{remark}

\bigskip {\noindent School of Mathematical Sciences, Tel-Aviv University, Tel-Aviv
69978, Israel \\  {\it e-mail address:}
\verb"roneneldan@gmail.com" }

\end{document}